%%%%%%%%%%%%%%%%%%%%%%%%%%%%%%%%%%%%%%%%%%%%%%%%%%%%%%%%%%%%%%%%%%%%%%%%%%
% This paper is being submitted to the Communications in Mathematical Physics.
% This version was completed on January 4, 1996.
%
% Title: Minimal model fusion rules from 2-groups
%
% Authors: Fusun Akman, Alex J. Feingold, Michael D. Weiner
% 
% Email Address: alex@math.binghamton.edu
%
% Phone Number: 607-777-2465  or  607-777-2147 (Secretary's Office)
% 
% Fax Number: 607-777-2450 
% 
% Mail Address: Professor Alex J. Feingold
%               Department of Mathematical Sciences
%               State University of New York
%               Binghamton, New York 13902-6000
%
% Subject Classifications: Primary: 17B68, 81R10
%                          Secondary: 81T40
%
% Since the second author was supported by a National Security Agency Grant,
% it is required that the following note be included with the submitted
% paper. 
% This manuscript is submitted for publication with the understanding that
% the United States Government is authorized to reproduce and distribute
% reprints.
%%%%%%%%%%%%%%%%%%%%%%%%%%%%%%%%%%%%%%%%%%%%%%%%%%%%%%%%%%%%%%%%%%%%%%%%%%
% This paper contains some tables which are produced by the TeX macro 
% cellular. The main file, cellular, and its related files were combined
% into a tar file by uufiles. After unpacking them, you should only have 
% to give the command to tex the main file, fusionrules.tex. I hope it works!
%%%%%%%%%%%%%%%%%%%%%%%%%%%%%%%%%%%%%%%%%%%%%%%%%%%%%%%%%%%%%%%%%%%%%%%%%%
\input amstex
\documentstyle{amsppt}
\magnification=\magstep1
\pageheight{9truein}
\pagewidth{6.5truein}
\NoRunningHeads
\NoBlackBoxes
\baselineskip=12pt
\def\f#1#2{{\tsize{\frac{#1}{#2}}}}

\def\bC{{\bold C}}
\def\boZ{{\bold Z}}
\def\cV{{\Cal V}}
\def\cN{{\Cal N}}
\def\Vir{\bold{Vir}}

\def\bk{\blacksquare}
\def\lra{\leftrightarrow}
\def\x{\times}

\let\pr\proclaim
\let\epr\endproclaim

\leftheadtext{AKMAN, FEINGOLD, WEINER}
\rightheadtext{MINIMAL MODEL FUSION RULES}

\topmatter
\title
Minimal Model Fusion Rules From $2$-Groups
\endtitle

\author 
F\"usun Akman \\ 
Alex J. Feingold 
\footnote{Partially supported by the National Security Agency under
Grant number MDA904-94-H-2019. The United States Government is authorized
to reproduce and distribute reprints notwithstanding any copyright notation
hereon. \hbox to 4.8truein{}
\endgraf
This paper is in final form and no version of it will be submitted for
publication elsewhere. \hbox to 1truein{}
\endgraf
Revised: Jan. 3, 1996 \hbox to 5.1truein{} 
}\\
Michael D. Weiner
\endauthor

\affil
Department of Mathematics \\
Cornell University\\
Ithaca, New York 14853-7901
\\ \\
Department of Mathematical Sciences \\
The State University of New York \\
Binghamton, New York 13902-6000
\\ \\
Department of Mathematics \\
Colgate University \\
Hamilton, NY 13346 
\endaffil

\email
alex\@math.binghamton.edu
\endemail


\endtopmatter

\document

\centerline{\bf{Abstract}}

The fusion rules for the $(p,q)$-minimal model representations of the Virasoro 
algebra are shown to come from the group $G = \boZ_2^{p+q-5}$ in the following
manner. There is a partition $G = P_1 \cup ...\cup P_N$ into disjoint subsets
and a bijection between $\{P_1,...,P_N\}$ and the sectors $\{S_1,...,S_N\}$
of the $(p,q)$-minimal model such that the fusion rules $S_i * S_j = \sum_k
D(S_i,S_j,S_k) S_k$ correspond to $P_i * P_j = \sum_{k\in T(i,j)} P_k$
where $T(i,j) = \{k|\exists a\in P_i,\exists b\in P_j, a+b\in P_k\}$.
\vskip 10pt

The Virasoro algebra $\Vir$ is the infinite dimensional Lie algebra with basis 
$\{c, L_m\ |\ m\in\boZ\}$ and brackets
$$[L_m,L_n] = (m-n)L_{m+n} + \frac{m^3-m}{12} \delta_{m,-n} c$$
where $c$ is central. A highest weight representation \cite{KR} of $\Vir$ is a
module $M$ containing a vector $v_0$ such that $L_m v_0 = 0$ for all $m\geq 1$,
$L_0 v_0 = hv_0$ for some $h\in{\bC}$, and $M$ is generated by $v_0$. 
It means that
$M$ is spanned by vectors of the form
$$L_{-n_1} L_{-n_2}\hdots L_{-n_r}v_0$$
for $n_1\geq n_2\geq \hdots \geq n_r\geq 1$ with $r\geq 0$. Since $c$ is
central, its value on $M$ is a scalar, which we continue to denote by $c$.
The Verma module $M(c,h)$ is the module, determined up to isomorphism by the
values of $c$ and $h$, such that the above vectors form a basis. For certain
values of these parameters the Verma module is not irreducible, but contains
nonzero null vectors, that is, vectors which generate proper highest weight
submodules. Taking the quotient of $M(c,h)$ by its maximal proper submodule,
consisting of the sum of all submodules generated by null vectors, yields the 
irreducible $\Vir$-module denoted by $V(c,h)$. Note that in order to
determine if a vector $v_0$ such that $L_0 v_0 = hv_0$ is null,
it suffices to check whether $L_1$ and
$L_2$ kill $v_0$, since the bracket formula shows that brackets of these
generate all $L_m$ for $m\geq 1$. 

In \cite{FZ} it is shown that for any $c\in{\bC}$, $V(c,0)$ has the structure 
of a Vertex Operator Algebra (VOA) \cite{FLM}
and each $V(c,h)$ is a module for that VOA. Various special   
classes of $\Vir$-modules have been investigated and classified. A highest 
weight module with $c$ and $h$ real admits a Hermitian form
$(\cdot,\cdot)$ determined by the conditions
$$(v_0,v_0) = 1\qquad\hbox{and}\qquad (L_m v,w) = (v,L_{-m}w)$$
for any $v,w$ in the module and any integer $m$. 
An important class of highest weight modules consists of the unitary modules 
for which the Hermitian form is positive definite. If
$c\geq 1$ and $h\geq 0$ are real then $V(c,h)$ is unitary, but there is also
a discrete series of unitary modules with $c < 1$ and only a finite number
of corresponding special $h$ values. The $c$ values are given by the formula
$$c = 1 - \frac{6}{p(p+1)} \hbox{ for } p\geq 2$$
and the $h$ values are given by
$$h_{m,n} = \frac{(n(p+1)-mp)^2 - 1}{4p(p+1)}\hbox{ for } 0<m<p+1,\ 0<n<p.$$
A larger discrete series of modules was investigated in \cite{BPZ}, and 
called the  minimal models. For integers $p,q\geq 2$ with $\gcd(p,q)=1$, the
$(p,q)$-minimal model consists of the Virasoro modules $V(c,h_{m,n})$ 
whose values of $c$ and $h_{m,n}$ are given by
$$c = 1 - \frac{6(p-q)^2}{pq}$$
and 
$$h_{m,n} = \frac{(np-mq)^2 - (p-q)^2}{4pq}\ \hbox{ for } 0<m<p,\ 0<n<q.$$
{}From the above formula, which we use to define $h_{m,n}$ for the rest of
the paper, we see that $h_{m,n} = h_{p-m,q-n}$, and 
the number of distinct such $h_{m,n}$ is $N = \f{1}{2} (p-1)(q-1)$.
These modules are of special interest from the viewpoint of VOAs because 
$V(c,0)$ has the structure of a rational VOA whose only irreducible 
VOA-modules are the modules $V(c,h_{m,n})$. 
In that case, one may hope to completely analyze the spaces of intertwining
operators between modules (see \cite{FHL}) and understand the fusion rules
which give the dimensions of those spaces. In \cite{W} the fusion rules for 
all minimal models are proved. The result is as follows.

\pr{Definition} For integer $p\geq 2$, we say that the triple of integers
$(m,m',m'')$ is {\bf p-admissible}
when  $0<m,m',m''<p$, the sum $m+m'+m''<2p$ is odd, and 
the ``triangle'' inequalities $m<m'+m''$, $m'<m+m''$, and $m''<m+m'$ are 
satisfied.
For $p,q\geq 2$ with $\gcd(p,q)=1$, we say that $((m,n),(m',n'),(m'',n''))$ 
is {\bf (p,q)-admissible} when $(m,m',m'')$ is $p$-admissible and 
$(n,n',n'')$ is $q$-admissible. \epr

\pr{Lemma 1} For $0<m'\leq m<p$, the set of all $m''$ such that $(m,m',m'')$ is
p-admissible equals $\{m-m'+1+2i\ |\ 0\leq i\leq \min(m',p-m)-1\}$. \epr 

\demo{Proof} If $0<m'\leq m<p$ then $m'' = m-m'+1+2i$ means $m+m'+m'' =
2m+1+2i$ is odd. Two of the triangle inequalities are satisfied: $m'+m''\geq
m'+(m-m'+1)=m+1>m$, and $m+m''\geq m+(m-m'+1)\geq m'+1>m'$.
We consider the two cases where (1) $m'\leq p-m$, and (2) $m'>p-m$. 

{\bf Case 1:} Suppose that $1\leq m-m'+1\leq m''\leq m+m'-1$ and $m+m'+m''$ is
odd. Then $m+m'>m+m'-1\geq m''$. We have $m''<p$ because
$m''\leq m+m'-1\leq m+(p-m)-1 =p-1<p$, so that $m+m'+m''\leq p+m''<2p$.

{\bf Case 2:} Let $1\leq m-m'+1\leq m''\leq 2p-(m+m')-1$ with $m+m'+m''$ odd.
We have $m''\leq 2p-(m+m')-1 < 2p-p-1 <p$ so $m+m'>p>m''$. We also have
$m+m'+m''\leq m+m'+2p-(m+m')-1 < 2p$.

We now have to show that $(m,m',m'')$ is not $p$-admissible if $m''<
m-m'+1$ or if $m'' > m-m'-1+2 \min(m',p-m)$. If $m''<m-m'+1$, then
$m'+m''< m+1$, so $m'+m''\leq m$, violating a triangle inequality. For large 
values of $m''$ we again consider the two cases above. In case (1) the condition
$m''>m+m'-1$ implies $m''\geq m+m'$. In case (2) if $m''>2p-(m+m')-1$ then 
$m+m'+m''> 2p-1$, so $m+m'+m''\geq 2p$. $\hfill\bk$

\enddemo

\pr{Definition} Let $\cN_{(m',n'),(m'',n'')}^{(m,n)} = 1$ if 
$((m,n),(m',n'),(m'',n''))$ is $(p,q)$-admissible, $0$ otherwise.\epr 

\pr{Theorem} [Wang] The fusion rules for the minimal models are
$$V(c,h_{m',n'}) \x V(c,h_{m'',n''}) = \sum_{(m,n)} 
\cN_{(m',n'),(m'',n'')}^{(m,n)} V(c,h_{m,n}).$$ 
\epr
Note that the definition of $p$-admissible is totally symmetric in the three
integers, so the definition of $(p,q)$-admissible has an $S_3\times S_3$
symmetry. While $h_{m,n} = h_{p-m,q-n}$, the definition of
admissibility does not allow the arbitrary replacement of a pair
$(m,n)$ by $(p-m,q-n)$. In fact, if the triple $((m,n),(m',n'),(m'',n''))$
is $(p,q)$-admissible, then $((m,n),(m',n'),(p-m'',q-n''))$ is not, as
$m+m'+m''$ and $n+n'+n''$ are odd, while $m+m'+p-m''$ and $n+n'+q-n''$
cannot both be odd ($p$, $q$ cannot have a common factor of 2). Therefore 
the sum in the Theorem is meant to be over all pairs $(m,n)$, but only
one of $V(c,h_{m,n})$ and $V(c,h_{p-m,q-n})$ may have a nonzero coefficient.

For a fixed $(p,q)$-minimal model, we will use the notation $[h] = V(c,h)$ 
to write the fusion rules more briefly. We may give these rules  
in a table, which may be interpreted as the
multiplication table for the associated Verlinde algebra. 
The fusion rules involving $[0]$ are just $[0]\x[h]=[h]$, but we
include these in our tables anyway. 

\pr{Definition} The Verlinde algebra $\cV$ associated with a fixed 
$(p,q)$-minimal model is the
associative algebra with basis indexed by the distinct modules $[h_{m,n}]$ 
and with product given by the fusion rules. \epr
 
In the case of $p=3$, $q=4$, one gets $c=\f{1}{2}$. 
The values of $h_{m,n} = h_{3-m,4-n}$ 
for $0<m<3$ and $0<n<4$ are given in the following table. 

\vskip 10pt
% File:       TeX Inputs cellular.tex
% Author:     J E Pittman
% Bitnet:     JEPTeX@TAMVenus
% Internet:   JEPTeX@Venus.TAMU.EDU
% Date:       November 8, 1988
%
% This file defines the main macro for cellular table construction.  
% For commentary, please see the file cellular.doc.
%
\message{cellular.TeX version 0.}%
\def\centertable{\leftskip=0pt plus1fill\rightskip=0pt plus1fill}

\def\begincellular#1#2\endcellular{\relax
   \begingroup
      % File:       TeX Inputs Cell1.tex
% Author:     J E Pittman
% Bitnet:     JEPTeX@TAMVenus
% Internet:   JEPTeX@Venus.TAMU.EDU
% Date:       October 11, 1988
%
% Set up the cellular environment
%
\catcode`_=11 % Protect local macros.
\ifx\forcount\undefined % File:       TeX Inputs loopy.tex
% Author:     J E Pittman
% Bitnet:     JEPTeX@TAMVenus
% Internet:   JEPTeX@Venus.TAMU.EDU
% Date:       September 29, 1988
%
% These macros supply structured, nested loop constructs as follows:
%  
%     \forcount \csname=number to number by number do
%        valid_TeX_input
%        \endfor \csname 
%
%     \while \csname \if_of_any_type do
%        valid_TeX_input
%        \endwhile \csname 
%
%     \whilenot \csname \if_of_any_type do
%        valid_TeX_input
%        \endwhile \csname         
%
% The \csname for the for loop must refer to a count register, for the 
% while loops it must uniquely identify the while loop.
%
\def\forcount #1{\relax
   \def
         \for #1=##1to ##2by ##3do
            ##4%
            \endfor #1%
      {\relax
      #1=##1\relax
      \ifnum ##3>0
         \whilenot #1\ifnum ##2<#1do
            ##4%
            \advance #1 by ##3\relax
            \endwhilenot #1%
      \else
         \while #1\ifnum ##2<#1do
            ##4%
            \advance #1 by ##3\relax
            \endwhile #1%
      \fi
      }%
   \for #1%   
   }%
\let\endwhilenot=\fi
\def\whilenot #1{\relax
   \def
         \whilenotloop#1 ##1do
            ##2%
            \endwhilenot #1%
      {\relax
         \expandafter\def\csname whilenotbody\string#1\endcsname{##2}%
         \expandafter\def\csname whilenotloop\string#1\endcsname
            {\relax
               ##1%
                  \let\next=\relax
               \else
                  \csname whilenotbody\string#1\endcsname
                  \expandafter\let\expandafter\next
                        \csname whilenotloop\string#1\endcsname
               \fi
               \next
               }%
         \csname whilenotloop\string#1\endcsname
         }%
   \whilenotloop#1
   }%
\let\endwhile=\fi
\def\while #1{\relax
   \def
         \whileloop#1 ##1do
            ##2%
            \endwhile #1%
      {\relax
         \expandafter\def\csname whilebody\string#1\endcsname{##2}%
         \expandafter\def\csname whileloop\string#1\endcsname
            {\relax
               ##1%
                  \csname whilebody\string#1\endcsname
                  \expandafter\let\expandafter\next
                        \csname whileloop\string#1\endcsname
               \else
                  \let\next=\relax
               \fi
               \next
               }%
         \csname whileloop\string#1\endcsname
         }%
   \whileloop#1
   }%

 \fi
\ifx\declarecount\undefined % File:       TeX Inputs declare.tex
% Author:     J E Pittman
% Bitnet:     JEPTeX@TAMVenus
% Internet:   JEPTeX@Venus.TAMU.EDU
% Date:       September 29, 1988
%
% These macros provide a method of locally allocating registers 
% without interference with previously allocated registers.  The 
% method is the same as on pages 346--347 of the \TeX book, however, 
% the declare macros are intended for local use only.  It is a logical 
% error to use a new macro between a declare macro and the end of the 
% appropriate enclosing group. 
%
\def\declarecount {\allocate0\countdef}%
\def\declaredimen {\allocate1\dimendef}%
\def\declareskip  {\allocate2\skipdef}%
\def\declaremuskip{\allocate3\muskipdef}%
\def\declarebox   {\allocate4\chardef}%
\def\declaretoks  {\allocate5\toksdef}%
\def\allocate#1#2#3{\relax
   \advance\count1#1 by 1
   \ifnum\count1#1<\count20
   \else
      \errmessage{No room for \string#3!}%
   \fi
   #2#3=\count1#1
   }%

 \fi
%
% Handy abbreviations
%
\def\half{0.5}%
\def\by{by}%
\def\height{height}%
\def\depth{depth}%
\def\width{width}%
\def\to{to}%
\def\zeropt{0pt}%
\def\notop{\toprulewidth=\zeropt\relax}%
\def\nobottom{\bottomrulewidth=\zeropt\relax}%
\def\noleft{\leftrulewidth=\zeropt\relax}%
\def\noright{\rightrulewidth=\zeropt\relax}%
\let\x_after=\expandafter
%
% When using the Xerox 9700s or 4050, use \setverticaladjustment for 
% portrait output and \sethorizontaladjustment for landscape output 
% due to the differences in the way that vertical and horizontal lines 
% of the same weight are printed.
%
\declaredimen\pixelwidth
\pixelwidth=1in
\divide\pixelwidth by 300                         % assume 300dpi
\declaredimen\horizontal_rule_adjust
\horizontal_rule_adjust=\zeropt
\def\sethorizontaladjustment{\horizontal_rule_adjust=\pixelwidth}%
\declaredimen\vertical_rule_adjust
\vertical_rule_adjust=\zeropt
\def\setverticaladjustment{\vertical_rule_adjust=\pixelwidth}%
%
% The left, right, bottom, and top rule widths are used to determine
% the widths of the box around each cell.
%
\declaredimen\leftrulewidth
\declaredimen\rightrulewidth
\declaredimen\bottomrulewidth
\declaredimen\toprulewidth
%
% The left, right, bottom, and top border skips are used to position 
% the text of a cell within it, relative to the centers of the rulers.
%
\declareskip\leftborderskip
\declareskip\rightborderskip
\declareskip\bottomborderskip
\declareskip\topborderskip
\declarecount\last_column
\declaredimen\columnwidth
\declarecount\merge_columns
\declaredimen\merge_width
\declarecount\last_row
\declaredimen\rowheight
\declarecount\merge_rows
\declaredimen\merge_height
\declarecount\rowpenalty
%
% The row info and column info token registers contain a list of 
% tokens of the form /number/info, where number is the number of a 
% row or column of interest and info is information, usually register 
% assignments, that pertains to the row or column.
%
\declaretoks\column_info
\column_info={/}%
\declaretoks\row_info
\row_info={/}%
\def\everycolumn{\leftrulewidth=0.4pt\relax
   \rightrulewidth=\leftrulewidth
   \leftborderskip=6pt plus 1fil\relax
   \rightborderskip=\leftborderskip
   \columnwidth=\zeropt\relax
   \merge_rows=0\relax
   \merge_height=\zeropt\relax
   \columnwidth=\zeropt\relax
   }%
\def\everyrow{\toprulewidth=0.4pt\relax
   \bottomrulewidth=\toprulewidth
   \topborderskip=3pt plus 1fil\relax
   \bottomborderskip=\topborderskip
   \rowheight=\zeropt\relax
   \merge_columns=0\relax
   \merge_width=\zeropt\relax
   }%
\def\get_data#1<#2{\relax
   \def\temp##1/#1/##2/##3***{\relax
      \def\temp{##2}%
      \ifnum1=0\temp
         #2={##1/#1//}%
      \else
%        \message{extracted ##2}% debug
         ##2%
      \fi
      }%
   \x_after\temp\the#2#1/1/***%
   }%
\def\add_data#1>#2#3{\relax
   \def\temp##1/#1/##2/##3***{\relax
      #2={##1/#1/##2#3/##3}%
%     \message{\string#2=\the#2}% debug
      }%
   \x_after\temp\the#2***%
   }%
\def\add_column_number_data{\relax
   \x_after \add_data \the\column_number>\column_info
   }%
\def\get_column_number_data{\relax
   \x_after \get_data \the\column_number<\column_info
   }%
\def\add_row_number_data{\relax
   \x_after \add_data \the\row_number>\row_info
   }%
\def\get_row_number_data{\relax
   \x_after \get_data \the\row_number<\row_info
   }%
\declarebox\temp_box
\declarebox\scratch_box
\declaredimen\temp_dimen
\declaredimen\scratch_dimen
\declareskip\temp_skip
\declarecount\temp_count
\declarecount\tracingexpansions
\tracingexpansions=0
\catcode`_=8 % Return to normal.
%
 
                                % set up enviroment
      #1\relax
      % File:       TeX Inputs Cell2.tex
% Author:     J E Pittman
% Bitnet:     JEPTeX@TAMVenus
% Internet:   JEPTeX@Venus.TAMU.EDU
% Date:       November 8, 1988
%
% Prepare to scan the data, taking notes as to span sizes, row and 
% column dimensions, et cetera.
%
\catcode`_=11 % Used to protect local control sequence names.
%
% The span info tokens contain sets of entries of the form \process 
% {position}{number}{dimension}, where position is the terminal column 
% or row, number is the number of columns or rows leading into the 
% column or row, and dimension is the size of the information.
% 
\declaretoks\column_span_info
\column_span_info={}%
\declaretoks\row_span_info
\row_span_info={}%
\let\process=\relax
\declarecount\column_number
\column_number=0
%
% Create a column information entry and put the user's specifications 
% into it.
%
\def\column#1{\relax
   \advance\column_number \by 1
   \last_column=\column_number
   \get_column_number_data
   \add_column_number_data {#1}%
   \ignorespaces
   }%
\declarecount\row_number
\row_number=0
%
% Same as \column.
%
\def\row#1{\relax
   \advance\row_number \by 1
   \message{Scanning row \the\row_number.}%
   \last_row=\row_number
   \everyrow
   \get_row_number_data
   \add_row_number_data {#1}%
   \column_number=0
   \ignorespaces
   }%
%
% \blank is used to generate a cell without a border or data.  In this 
% context, all it does is absorb merges.
%
\def\blank{\relax
   \advance\column_number \by 1
   \if\column_number>\last_column
      \advance\column_number \by -1
      \column{}%
   \fi
   \everycolumn
   \get_column_number_data
   \ifnum\merge_rows>1
      \add_column_number_data {\merge_rows=0\relax}%
   \fi
   \merge_columns=0
   }%
%
% \cell is used to generate a normal, ruled cell.  In this context, it 
% merely measures the cell and makes the appropriate notes.
%
\def\cell#1{\relax
   \advance\column_number \by 1
   \if\column_number>\last_column
      \advance\column_number \by -1
      \column{}%
   \fi
   \everycolumn
   \get_column_number_data
%
% Typeset the information into temp box.
%
   \setbox\temp_box=\vbox \bgroup
      \begingroup
         \ifnum\merge_rows>0
            \advance\row_number \by -\merge_rows
            \get_row_number_data
         \fi
         \vskip \topborderskip
         \endgroup
      \hbox \bgroup
         \begingroup
            \ifnum\merge_columns>0
               \advance\column_number \by -\merge_columns
               \get_column_number_data
            \fi
            \hskip \leftborderskip
            \endgroup
         #1\vphantom{)}%
         \hskip \rightborderskip
         \egroup
      \vskip \bottomborderskip
      \egroup
%
% If it is a row merger, record it for later processing.
%
   \ifnum \merge_rows>0
      \edef\temp{\process
         {\the\merge_rows}{\the\ht\temp_box}{\the\row_number}%
         \the\row_span_info
         }%
      \x_after\row_span_info\x_after=\x_after{\temp}%
%     \message{\string\row_span_info=\the\row_span_info}% debug
      \add_column_number_data {\merge_rows=0\relax}%
   \else
%
% Not a merger, record the height if max.
%
      \ifdim\ht\temp_box>\rowheight
         \let\info=\relax
         \edef\temp{\the\row_number>\info
               {\rowheight=\the\ht\temp_box\relax}}%
         \let\info=\row_info
         \x_after \add_data \temp
%        \message{\string\row_info=\the\row_info}% debug
         \rowheight=\ht\temp_box
      \fi
   \fi
%
% Same as above for column merger and width.
%
   \ifnum \merge_columns>0
      \edef\temp{\process
         {\the\merge_columns}{\the\wd\temp_box}{\the\column_number}%
         \the\column_span_info
         }%
      \x_after\column_span_info\x_after=\x_after{\temp}%
%     \message{\string\column_span_info=\the\column_span_info}% debug
      \merge_columns=0
   \else
      \ifdim\wd\temp_box>\columnwidth
         \let\info=\relax
         \edef\temp{\the\column_number>\info
               {\columnwidth=\the\wd\temp_box\relax}}%
         \let\info=\column_info
         \x_after \add_data \temp
%        \message{\string\column_info=\the\column_info}% debug
      \fi
   \fi
   }%
%
% dp 2/19 -- added the following definitions for left justified, centered,
% and right justified cells
%
\def\rcell#1{\cell{\hfill{}#1}}
\def\ccell#1{\cell{\hfill{}#1\hfill{}}}
\def\lcell#1{\cell{#1\hfill{}}}

%
% \mergeright specifies that the corresponding position is to be 
% merged with the cell to its right.
%
\def\mergeright{\relax
   \advance\column_number \by 1
   \if\column_number>\last_column
      \advance\column_number \by -1
      \column{}%
   \fi
   \everycolumn
   \get_column_number_data
   \advance\merge_columns \by 1
%
% Cancel a row merge, if present.
%
   \ifnum\merge_rows>1
      \add_column_number_data {\merge_rows=0\relax}%
   \fi
   }%
%
% Same as \mergeright, except down.
%
\def\mergedown{\relax
   \advance\column_number \by 1
   \if\column_number>\last_column
      \advance\column_number \by -1
      \column{}%
   \fi
   \everycolumn
   \get_column_number_data
   \add_column_number_data {\advance\merge_rows \by 1\relax}%
   \merge_columns=0
   }%
%
% The horizontal and vertical stretch macros allow the user to specify 
% an explicit stretch that will subsequently be processed like a span.  
% User-specified stretches are processed after span caused ones.  The 
% parameters are the starting column/row, the ending column/row, and 
% the size of the stretch.
%
\def\horizontalstretch#1#2#3{\relax
   \temp_count=#2\relax
   \advance\temp_count \by -#1\relax
   \edef\temp{\the\column_span_info\process{\the\temp_count}{#3}{#2}}%
   \x_after \column_span_info\x_after=\x_after{\temp}%
%  \message{\string\column_span_info=\the\column_span_info}% debug
   \ignorespaces
   }%
\def\verticalstretch#1#2#3{\relax
   \temp_count=#2\relax
   \advance\temp_count \by -#1\relax
   \edef\temp{\the\row_span_info\process{\the\temp_count}{#3}{#2}}%
   \x_after \row_span_info\x_after=\x_after{\temp}%
%  \message{\string\row_span_info=\the\row_span_info}% debug
   \ignorespaces
   }%
\def\noalign#1{\ignorespaces}% don't do anything for the first pass
\catcode`_=8 % Return to normal.
%
 
                                % set up for scan
      \ignorespaces
      #2\relax                                    % scan the cells
      % File:       TeX Inputs Cell3.tex
% Author:     J E Pittman
% Bitnet:     JEPTeX@TAMVenus
% Internet:   JEPTeX@Venus.TAMU.EDU
% Date:       October 11, 1988
%
% Process the column and row span info.
%
\declaredimen\expansion
\edef\everycolumn{\everycolumn\expansion=\zeropt\relax}%
\edef\everyrow{\everyrow\expansion=\zeropt\relax}%
\catcode`_=11 % used to protect local control sequence names.
%
%\message{\string\row_span_info=\the\row_span_info}% debug
%\message{\string\column_span_info=\the\column_span_info}% debug
%
\def\process#1#2#3{\relax
   \last_cell=#3\relax
   \first_cell=\last_cell
   \advance \first_cell \by -#1\relax
   \span_size=#2\relax
%
% Compute the gap between the size of the span and the total size of 
% the cells spanned.
%
   \gap=\span_size
   \forcount \cell_number=\first_cell to \last_cell by 1 do
      \everycell
      \get_cell_number_data
      \advance \gap \by -\cell_size
      \advance \gap \by -\expansion
      \endfor \cell_number
%  \message{\string\first_cell=\the\first_cell}% debug
%  \message{\string\last_cell=\the\last_cell}% debug
%  \message{\string\span_size=\the\span_size}% debug
%  \message{\string\gap=\the\gap}% debug
%
% If the gap is 0pt or less, nothing needs to be done, else search for 
% the minimum expansion that can be applied to every cell with a 
% current expansion less than the expansion found such that the span 
% is properly accomadated.
%
   \ifdim \gap>\zeropt
      \expandable_cells=#1\relax
      \advance \expandable_cells \by 1
      \trial_expansion=\zeropt
      \whilenot\search \ifdim\gap=\zeropt do
%        \message{\string\gap=\the\gap}% debug
         \ifnum \expandable_cells=0
            \advance \trial_expansion \by \expansion
         \else
            \multiply \trial_expansion \by \expandable_cells
            \advance \trial_expansion \by \gap
            \divide \trial_expansion \by \expandable_cells
            \expandable_cells=0
         \fi
         \gap=\span_size
%        \message{\string\trial_expansion=\the\trial_expansion}% debug
         \forcount \cell_number=\first_cell to \last_cell by 1 do
            \everycell
            \get_cell_number_data
            \advance \gap \by -\cell_size
            \ifdim \expansion>\trial_expansion
               \advance \gap \by -\expansion
            \else
               \advance \gap \by -\trial_expansion
               \advance \expandable_cells \by 1
            \fi
            \endfor \cell_number
         \temp_dimen=1sp
         \multiply \temp_dimen \by \expandable_cells
         \ifdim \gap>-\temp_dimen
            \ifdim \gap<\temp_dimen
               \gap=\zeropt
            \fi
         \fi
         \endwhilenot \search
      \forcount \cell_number=\first_cell to \last_cell by 1 do
         \everycell
         \get_cell_number_data
         \ifdim \expansion<\trial_expansion
            \let\info=\relax
            \edef\temp{\the\cell_number>\info
                  {\expansion=\the\trial_expansion\relax}}%
            \let\info=\cell_info
            \x_after \add_data \temp
            \ifnum\tracingexpansions>0
               \message{Expanded \the\cell_number}%
               \message{by \the\trial_expansion}%
               \message{from \the\cell_size}%
               \advance \cell_size \by \trial_expansion
               \message{to \the\cell_size.}%
            \fi
         \fi
         \endfor \cell_number
   \fi
   }%
\declarecount\first_cell
\declarecount\last_cell
\declaredimen\span_size
\let\expandable_cells=\temp_count
\declaredimen\trial_expansion
\let\gap=\scratch_dimen
\let\cell_number=\row_number
\let\everycell=\everyrow
\let\get_cell_number_data=\get_row_number_data
\let\cell_info=\row_info
\let\cell_size=\rowheight
\ifnum\tracingexpansions>0
   \message{Checking row expansions.}%
\fi
\the\row_span_info
\let\cell_number=\column_number
\let\everycell=\everycolumn
\let\get_cell_number_data=\get_column_number_data
\let\cell_info=\column_info
\let\cell_size=\columnwidth
\ifnum\tracingexpansions>0
   \message{Checking column expansions.}%
\fi
\the\column_span_info
\let\process=\relax
\catcode`_=8 % back to normal
%
 
                                % compute spans
      % File:       TeX Inputs Cell4.tex
% Author:     J E Pittman
% Bitnet:     JEPTeX@TAMVenus
% Internet:   JEPTeX@Venus.TAMU.EDU
% Date:       November 8, 1988
%
% Set up to output the data.
%
\catcode`_=11 % Protect local control sequence names.
%
% The user supplied information about the column has already been 
% processed.
%
\def\column #1{\relax\ignorespaces}%
\row_number=0
\rowpenalty=0
%
% This routine is used for horizontal kerning when there might be a 
% kern to the left of the current position.
%
\def\move_right_via_lastkern #1{\relax
   \temp_dimen=#1\relax
   \ifdim \lastkern>\zeropt
      \advance \temp_dimen \by \lastkern
      \unkern
   \else
   \fi
   \kern \temp_dimen
   }%
%
% \row begins a row by getting its specifications, terminating the 
% previous row (if any) and going into horizontal mode.
%
\def\row #1{\relax
   \advance \row_number \by 1
   \everyrow
   \get_row_number_data
   \advance \rowheight \by \expansion
   \ifdim \bottomrulewidth>\zeropt
      \advance \bottomrulewidth \by \horizontal_rule_adjust
   \fi
   \column_number=0
   \par
   \ifnum \rowpenalty=0
   \else
      \penalty \rowpenalty
      \rowpenalty=0
   \fi
   \noindent
   \ignorespaces
   \message{Outputting row \the\row_number.}%
   }%
%
% \blank creates a blank cell by kerning the appropriate amount.
%
\def\blank {\relax
   \advance \column_number \by 1
   \everycolumn
   \get_column_number_data
   \advance \columnwidth \by \expansion
   \advance \merge_width \by \expansion
   \move_right_via_lastkern \merge_width
%
% Terminate merger(s).
%
   \merge_width=\zeropt
   \merge_columns=0
   \ifnum \merge_rows>0
      \add_column_number_data
            {\merge_rows=0\relax\merge_height=\zeropt\relax}%
   \fi
   }%
%
% \cell outputs a cell.  The components of the cell are (in the order 
% output) the entry, the top ruler, the bottom ruler, and the left and 
% right rulers.
%
\def\cell #1{\relax
   \advance \column_number \by 1
   \everycolumn
   \get_column_number_data
   \advance \columnwidth \by \expansion
   \advance \merge_height \by \rowheight
   \advance \merge_width \by \columnwidth
   \ifdim \leftrulewidth>\zeropt
      \advance \leftrulewidth \by \vertical_rule_adjust
   \fi
   \ifdim \rightrulewidth>\zeropt
      \advance \rightrulewidth \by \vertical_rule_adjust
   \fi
%
% Get the correct top border skip and rule width.  Note that it is 
% necessary to extract this informaion even if a row merger is not 
% present because a previous row merger might have left the wrong 
% values.
%
   \begingroup
      \advance \row_number \by -\merge_rows
      \everyrow
      \get_row_number_data
      \xdef\globaltemp{\topborderskip=\the\topborderskip\relax
         \toprulewidth=\the\toprulewidth\relax
         }%
      \aftergroup \globaltemp
      \endgroup
   \ifdim \toprulewidth>\zeropt
      \advance \toprulewidth \by \horizontal_rule_adjust
   \fi
%
% Same procedure for the left border skip and rule width except that 
% extraction is necessary only in the presense of a column merger due 
% to the execution of an every column and a get at the start of \cell.  
%
   \ifnum \merge_columns>0
      \begingroup
         \advance \column_number \by -\merge_columns
         \everycolumn
         \get_column_number_data
         \xdef\globaltemp{\leftrulewidth=\the\leftrulewidth\relax
            \leftborderskip=\the\leftborderskip\relax
            }%
         \aftergroup \globaltemp
         \endgroup
      \ifdim \leftrulewidth>\zeropt
         \advance \leftrulewidth \by \vertical_rule_adjust
      \fi
   \fi
%
% Typeset the entry into temp box horizontally first, trying kerns 
% before glue in case the cell does not require horizontal stretching 
% and taking advantage of an empty cell by doing nothing, if such is 
% the case.
%
   \setbox\temp_box=\hbox{#1}%
   \ifdim\wd\temp_box>\zeropt
      \setbox\temp_box=\hbox \bgroup
         \kern \leftborderskip
         \box\temp_box
         \egroup
      \temp_dimen=\wd\temp_box
      \advance\temp_dimen \by \rightborderskip
      \wd\temp_box=\temp_dimen
      \ifdim\wd\temp_box=\merge_width
%
% then the kerns can be used instead of skips.
%
      \else
         \setbox\temp_box=\hbox \to \merge_width \bgroup
            \hskip \leftborderskip
            #1%
            \hskip \rightborderskip
            \egroup
      \fi
%
% Hide the width of temp box and put a phantom into it the hard way.
%
      \wd\temp_box=\zeropt
      \setbox\scratch_box=\hbox{#1)}%
      \ifdim \dp\scratch_box>\dp\temp_box
         \dp\temp_box=\dp\scratch_box
      \fi
      \ifdim \ht\scratch_box>\ht\temp_box
         \ht\temp_box=\ht\scratch_box
      \fi
      \temp_dimen=\ht\temp_box
      \advance \temp_dimen \by \dp\temp_box
      \advance \temp_dimen \by \bottomborderskip
      \advance \temp_dimen \by \topborderskip
      \ifdim \temp_dimen=\merge_height
%
% then the entry can be positioned vertically via a raise statement.  
% The total height of the material output should be equal to the row 
% height, thus acting as a strut.
%
         \temp_dimen=\bottomborderskip
         \advance \temp_dimen \by \dp\temp_box
         \scratch_dimen=\rowheight
         \advance\scratch_dimen by -\temp_dimen
         \ht\temp_box=\scratch_dimen
         \raise \temp_dimen \box\temp_box
      \else % have to do it via a box
         \setbox\temp_box=\vbox \to \rowheight \bgroup
%
% Subtracting merge height - row height from top border skip allows 
% the cell to stick up into the next row by an appropriate amount.
%
            \advance \topborderskip \by \rowheight
            \advance \topborderskip \by -\merge_height
            \vskip \topborderskip
            \box\temp_box
            \vskip \bottomborderskip
            \egroup
         \box\temp_box
      \fi
   \fi
%
% All of the rules are typeset with an overlap of at least pixel width 
% which insures that there will be no gaps.
%
% Typeset the top rule into an hbox and use a raise statement to put 
% it into position.
%
   \ifdim \toprulewidth>\zeropt
      \setbox\temp_box=\hbox \bgroup
         \temp_dimen=\merge_width
         \ifdim \half\leftrulewidth<\pixelwidth
            \kern -\pixelwidth
         \else
            \kern -\half\leftrulewidth
         \fi
         \advance \temp_dimen \by -\lastkern
         \vrule \height \half\toprulewidth
                \depth  \half\toprulewidth
                \width  \temp_dimen
         \ifdim \half\rightrulewidth<\pixelwidth
            \temp_dimen=\pixelwidth
         \else
            \temp_dimen=\half\rightrulewidth
         \fi
         \kern -\temp_dimen
         \vrule \height \half\toprulewidth
                \depth  \half\toprulewidth
                \width  2\temp_dimen
         \egroup
      \wd\temp_box=\zeropt
      \temp_dimen=\rowheight
      \advance\temp_dimen \by -\merge_height
      \ht\temp_box=\temp_dimen
      \dp\temp_box=\merge_height
      \raise \merge_height \box\temp_box
   \fi
%
% Output the bottom rule using the same methods.
%
   \ifdim \bottomrulewidth>\zeropt
      \setbox\temp_box=\hbox \bgroup
         \temp_dimen=\merge_width
         \ifdim \half\leftrulewidth<\pixelwidth
            \kern -\pixelwidth
         \else
            \kern -\half\leftrulewidth
         \fi
         \advance \temp_dimen \by -\lastkern
         \vrule \height \half\bottomrulewidth
                \depth  \half\bottomrulewidth
                \width  \temp_dimen
         \ifdim \half\rightrulewidth<\pixelwidth
            \temp_dimen=\pixelwidth
         \else
            \temp_dimen=\half\rightrulewidth
         \fi
         \kern -\temp_dimen
         \vrule \height \half\bottomrulewidth
                \depth  \half\bottomrulewidth
                \width  2\temp_dimen
         \egroup
      \wd\temp_box=\zeropt
      \dp\temp_box=\zeropt
      \ht\temp_box=\rowheight
      \box\temp_box
   \fi
%
% Test to see if the left inclusive-or right rule width is non-zero.
%
   \ifdim \leftrulewidth=\zeropt
      \temp_dimen=\rightrulewidth
   \else
      \temp_dimen=\leftrulewidth
   \fi
   \ifdim \temp_dimen>\zeropt
      \setbox\temp_box=\hbox \bgroup
         \temp_dimen=\merge_height
         \advance \merge_height \by \pixelwidth
         \ifdim \leftrulewidth>\zeropt
            \kern -\half\leftrulewidth
            \vrule \height \temp_dimen
                   \depth  \pixelwidth
                   \width  \leftrulewidth
         \fi
         \ifdim \rightrulewidth>\zeropt
            \scratch_dimen=\merge_width
            \advance \scratch_dimen \by -\half\leftrulewidth
            \advance \scratch_dimen \by -\half\rightrulewidth
            \kern \scratch_dimen
            \vrule \height \temp_dimen
                   \depth  \pixelwidth
                   \width  \rightrulewidth
         \fi
         \egroup
      \wd\temp_box=\merge_width
      \ht\temp_box=\rowheight
      \dp\temp_box=\zeropt
      \box\temp_box
   \else
      \move_right_via_lastkern \merge_width
   \fi
%
% Cancel the mergers.
%
   \merge_width=\zeropt
   \merge_columns=0
   \ifnum \merge_rows>0
      \add_column_number_data
            {\merge_rows=0\relax\merge_height=\zeropt\relax}%
   \fi
   \ignorespaces
   }%
%
% No surprises here.
%
\def\mergeright {\relax
   \advance \column_number \by 1
   \everycolumn
   \get_column_number_data
   \advance \columnwidth \by \expansion
   \advance \merge_width \by \columnwidth
   \advance \merge_columns \by 1
   \ifnum \merge_rows>0
      \add_column_number_data
            {\merge_rows=0\relax\merge_height=\zeropt\relax}%
   \fi
   }%
%
% No surprises here.
%
\def\mergedown {\relax
   \advance \column_number \by 1
   \everycolumn
   \get_column_number_data
   \advance \columnwidth \by \expansion
   \advance \merge_width \by \columnwidth
   \move_right_via_lastkern \merge_width
   \merge_width=\zeropt
   \merge_columns=0
   \advance \merge_height \by \rowheight
   \let\info=\relax
   \edef\temp{\the\column_number>\info
         {\merge_height=\the\merge_height\relax
         \advance\merge_rows \by 1\relax}}%
   \let\info=\column_info
   \x_after \add_data \temp
   \rowpenalty=10000 % do not allow a break over a row merge.
   }%
\catcode`_=8 % Back to normal.
\def\noalign#1{\relax
   \vadjust{#1}%
   \ignorespaces
   }%
%
 
                                % set up for output
      \offinterlineskip
      \parskip=\zeropt
      \ignorespaces
      #2\relax                                    % output cells
      \par
      \endgroup
   }%

\begincellular{\centertable}
\row{}\cell{$h_{m,n}$}\cell{$n=1$}\cell{$n=2$}\cell{$n=3$}
\row{}\cell{$m=1$}\cell{$0$}\cell{$\f{1}{16}$}\cell{$\f{1}{2}$}
\row{}\cell{$m=2$}\cell{$\f{1}{2}$}\cell{$\f{1}{16}$}\cell{$0$}
\endcellular
\vskip 15pt

The $c = \f{1}{2}$ fusion rule table is as follows:

\vskip 10pt
\begincellular{\centertable}
\row{}
\cell{$[h]\x[h']$}\cell{$[0]$}\cell{$[\f{1}{16}]$}\cell{$[\f{1}{2}]$}
\row{}
\cell{$[0]$}
\cell{$[0]$}
\cell{$[\f{1}{16}]$}
\cell{$[\f{1}{2}]$}
\row{}
\cell{$[\f{1}{16}]$}
\cell{$[\f{1}{16}]$}
\cell{$[0]+[\f{1}{2}]$}
\cell{$[\f{1}{16}]$}
\row{}
\cell{$[\f{1}{2}]$}
\cell{$[\f{1}{2}]$}
\cell{$[\f{1}{16}]$}
\cell{$[0]$}
\endcellular
\vskip 15pt

This well-known case, called the Ising model, was studied in great detail in 
\cite{FRW} using the spinor construction \cite{FFR} from one fermion. 
The objective was to unify the  
superVOA $V(\f{1}{2},0)\oplus V(\f{1}{2},\f{1}{2})$ 
(the Neveu-Schwarz sector)
with its module $V(\f{1}{2},\f{1}{16})\oplus V(\f{1}{2},\f{1}{16})$ 
(the Ramond sector), and to
find the generalization of the Jacobi-Cauchy identity obeyed by the
intertwining operators from the Ramond sector. The spinor construction
naturally provided two copies of the $h=\f{1}{16}$ module, and the above Ising
fusion rules were replaced by fusion rules given by the group ${\boZ}_4$. It
meant that each vector in the Ramond sector corresponded to a unique
intertwining operator. The question remained whether that was a special
situation, or if other fusion rules could be rewritten in terms of a group
by taking multiple copies of certain sectors. While the Ising fusion rules
could also have come from the group
${\boZ}_2\times {\boZ}_2$, it was not checked if that choice would be
consistent with the VOA-module structure and the generalization of the
Jacobi-Cauchy identity for intertwiners. 

Let us explain further what we mean by the replacement of these Ising fusion
rules by fusion rules given by the group ${\boZ}_4 = \{0,1,2,3\}$ under
addition modulo 4. We associate these four group elements with the four
irreducible $\Vir$-modules from the fermionic construction as follows: 
$$0\lra [0],\quad 1\lra [\f{1}{16}]_1,\quad 2\lra [\f{1}{2}],
\quad 3\lra [\f{1}{16}]_2.$$
Then the addition in ${\boZ}_4$ would give the Ising fusion rules if we
could not distinguish between the two copies of the $[\f{1}{16}]$-module.
Below we will give a precise definition of when a finite abelian group
covers the fusion rules of a $(p,q)$-minimal model.

The next example to examine is $p=4$, $q=5$, so $c=\f{7}{10}$,
$h\in\{0,\f{3}{5},\f{7}{16},\f{3}{80},\f{1}{10},\f{3}{2}\}$. 
The values of $h_{m,n} = h_{4-m,5-n}$ 
for $0<m<4$ and $0<n<5$ are given in the following table. 

\vskip 10pt

\begincellular{\centertable}
\row{}\cell{$h_{m,n}$}\cell{$n=1$}\cell{$n=2$}\cell{$n=3$}\cell{$n=4$}
\row{}\cell{$m=1$}\cell{$0$}\cell{$\f{1}{10}$}\cell{$\f{3}{5}$}
\cell{$\f{3}{2}$}
\row{}\cell{$m=2$}\cell{$\f{7}{16}$}\cell{$\f{3}{80}$}\cell{$\f{3}{80}$}
\cell{$\f{7}{16}$}
\row{}\cell{$m=3$}\cell{$\f{3}{2}$}\cell{$\f{3}{5}$}\cell{$\f{1}{10}$}
\cell{$0$}
\endcellular
\vskip 15pt

The $c = \f{7}{10}$ fusion rule table is as follows:

\vskip 10pt
\begincellular{\centertable}
\row{}
\cell{$[h]\x[h']$}\cell{$[0]$}\cell{$[\f{3}{5}]$}\cell{$[\f{7}{16}]$}
\cell{$[\f{3}{80}]$}\cell{$[\f{1}{10}]$}\cell{$[\f{3}{2}]$}
\row{}
\cell{$[0]$}
\cell{$[0]$}
\cell{$[\f{3}{5}]$}
\cell{$[\f{7}{16}]$}
\cell{$[\f{3}{80}]$}
\cell{$[\f{1}{10}]$}
\cell{$[\f{3}{2}]$}
\row{}
\cell{$[\f{3}{5}]$}
\cell{$[\f{3}{5}]$}
\cell{$[0]+[\f{3}{5}]$}
\cell{$[\f{3}{80}]$}
\cell{$[\f{3}{80}]+[\f{7}{16}]$}
\cell{$[\f{1}{10}]+[\f{3}{2}]$}
\cell{$[\f{1}{10}]$}
\row{}
\cell{$[\f{7}{16}]$}
\cell{$[\f{7}{16}]$}
\cell{$[\f{3}{80}]$}
\cell{$[0]+[\f{3}{2}]$}
\cell{$[\f{1}{10}]+[\f{3}{5}]$}
\cell{$[\f{3}{80}]$}
\cell{$[\f{7}{16}]$}
\row{}
\cell{$[\f{3}{80}]$}
\cell{$[\f{3}{80}]$}
\cell{$[\f{3}{80}]+[\f{7}{16}]$}
\cell{$[\f{1}{10}]+[\f{3}{5}]$}
\cell{$[0]+[\f{3}{5}]+[\f{3}{2}]+[\f{1}{10}]$}
\cell{$[\f{7}{16}]+[\f{3}{80}]$}
\cell{$[\f{3}{80}]$}
\row{}
\cell{$[\f{1}{10}]$}
\cell{$[\f{1}{10}]$}
\cell{$[\f{1}{10}]+[\f{3}{2}]$}
\cell{$[\f{3}{80}]$}
\cell{$[\f{7}{16}]+[\f{3}{80}]$}
\cell{$[0]+[\f{3}{5}]$}
\cell{$[\f{3}{5}]$}
\row{}
\cell{$[\f{3}{2}]$}
\cell{$[\f{3}{2}]$}
\cell{$[\f{1}{10}]$}
\cell{$[\f{7}{16}]$}
\cell{$[\f{3}{80}]$}
\cell{$[\f{3}{5}]$}
\cell{$[0]$}
\endcellular
\vskip 15pt

We would like to point out an interesting feature of the
above fusion rule table. If you examine the number of times that an entry
appears on each row, you find that $[0]$ and $[\f{3}{2}]$ each appear at most 
once, $[\f{3}{5}]$, $[\f{7}{16}]$ and $[\f{1}{10}]$ each appear at most twice, 
and $[\f{3}{80}]$ appears at most four times. In the cyclic group $\boZ_{12}$ 
we find a matching list of repetitions of elements of orders dividing $12$. 
Only 0 has order 1 and only 6 has order 2, but 4 and 8 have order 3, 3 and 9
have order 4, 2 and 10 have order 6, and 1, 5, 7 and 11 each have order 12.
If we take multiple copies of the modules, the multiplicity given by the
number of repetitions, we can get a perfect correspondence with the group
$\boZ_{12}$ so that the fusion rules are replaced by the group law of
addition modulo $12$. The correspondence is
$$\eqalign{
&0\lra [0],\ 1\lra [\f{3}{80}]_1,\ 2\lra [\f{1}{10}]_1, 3\lra [\f{7}{16}]_1,\
4\lra [\f{3}{5}]_1,\ 5\lra [\f{3}{80}]_2,\cr
&6\lra [\f{3}{2}], 7\lra [\f{3}{80}]_3,\ 
8\lra [\f{3}{5}]_2,\ 9\lra [\f{7}{16}]_2,\ 10\lra [\f{1}{10}]_2, 
11\lra [\f{3}{80}]_4.\cr}$$
We found that while some other minimal model fusion rules could be similarly
covered by finite cyclic groups, in order to cover all of them we had to use
certain $2$-groups. We suspect that this may be explained by the
Goddard-Kent-Olive coset construction \cite{GO}. 
This investigation is a first step toward a uniform
method of constructing all the minimal models, including the intertwining
operators, in such a way that each intertwining operator corresponds to a
unique vector in some module. The fusion rules will be replaced by a group
law, so the Verlinde algebra will be replaced by a group algebra. (This does
not mean there is an algebra map from the new group algebra to the old
Verlinde algebra.) We believe this will make possible a construction of an
algebraic system containing the direct sum of $V(c,0)$ and copies of each of 
its modules $V(c,h)$, so that the intertwiners are in one-to-one correspondence
with the states in the modules, and they obey a matrix generalization of the
Jacobi-Cauchy identity. 

Fix $p,q \ge 2$ with $\gcd(p,q)=1$. Let 
$$S(p,q) = \{(m,n)\ |\ 0 < m < p, 0 < n < q\}$$ 
and let the set of distinct modules in the $(p,q)$-minimal model 
be denoted by
$$V = V(p,q) = \{[h_{m,n}]\ |\ (m,n)\in S(p,q)\}.$$

\pr{Definition} Let $G$ be a finite abelian group. We say a map $\Phi : G\to
V(p,q)$ covers the fusion rules of the $(p,q)$-minimal model if the
following conditions are satisfied:

\item{(1)} If $g_i\in G$ and $\Phi(g_i) = [h_{m_i,n_i}]$ for $1\leq i\leq 2$,
then $\Phi(g_1 + g_2) = [h_{m_3,n_3}]$ for some $(m_3,n_3)\in S(p,q)$ with
$((m_1,n_1),(m_2,n_2),(m_3,n_3))$ $(p,q)$-admissible.

\item{(2)} If $((m_1,n_1),(m_2,n_2),(m_3,n_3))$ is $(p,q)$-admissible then
for $1\leq i\leq 3$ there exist elements $g_i\in G$ such that
$\Phi(g_i) = [h_{m_i,n_i}]$ and $g_1 + g_2 = g_3$.
\epr

This definition can also be described in the following alternative way.
Let $N$ be the number of distinct irreducible Virasoro
modules in the $(p,q)$-minimal model. 
Let $S_1,...,S_N$ be labels for these $N$ sectors, with $S_1$ corresponding to
$h_{1,1} = 0$, and let the fusion rules be denoted by $D(S_i,S_j,S_k)$. 
Then the Verlinde algebra $\cV$ is an $N$-dimensional vector space
over {\bf Q} with basis $\{S_1,...,S_N\}$ and with 
product $S_i * S_j = \sum_k D(S_i,S_j,S_k) S_k$.
Let $(G,+,0)$ be a finite abelian group and let
$G = P_1 \cup ... \cup P_N$ be a partition into 
$N$ disjoint subsets with $P_1 = \{0\}$.
Let $W$ be an $N$-dimensional vector space over {\bf Q} with basis
$\{P_1,...,P_N\}$ and give $W$ an algebra structure by defining the
product $P_i * P_j = \sum_{k\in T(i,j)} P_k$
where $T(i,j) = \{k|\exists a\in P_i,\exists b\in P_j, a+b\in P_k\}$.
We say that this partition of $G$ gives the fusion rules of the $(p,q)$-minimal
model if the bijection $S_i \leftrightarrow P_i$ determines an algebra 
isomorphism between $\cV$ and $W$.

Let $r = p+q-4$, and define the abelian group $H$ to be ${\boZ}_2^r$.
For $x = (x_1,\cdots,x_r)\in H$ define $supp(x)$, the  
{\it{support}} of $x$, to be the set of coordinates of $x$ which are $1$,
and define $wt(x) = |supp(x)|$, the {\it{weight}} of $x$, 
to be the number of coordinates of $x$ which are $1$. 
Of course, $H$ also has a Boolean ring structure with the product given by
coordinate-wise multiplication.

Define the subgroups $A = \{x\in H\ |\ supp(x)\subseteq \{1,\hdots,p-2\}\}$ and
$B = \{x\in H\ |\ supp(x)\subseteq \{p-1,\hdots,p+q-4\}\}$ in $H$.
Note that in $H$, $x+y = z$ is equivalent to $x+y+z = 0$, which has an obvious
$S_3$ symmetry. Also, if $a_i\in A$, $b_i\in B$, $1\leq i\leq 3$, then 
$(a_1+b_1) + (a_2+b_2) = a_3+b_3$ iff $(a_1+b_1) + (a_2+b_2) + (a_3+b_3) = 0$
iff $a_1 + a_2 + a_3 = 0$ and $b_1 + b_2 + b_3 = 0$, which has an 
$S_3\times S_3$ symmetry. The symmetric group $S_{p-2}$ acts as a
group of automorphisms of $A$ by permuting the coordinates, as $S_{q-2}$ acts 
on $B$. For $0 < m < p$ and $0 < n < q$, let
$$A_m = \{x \in A\ |\ wt(x) = m-1\}\quad\hbox{ and }\quad
B_n = \{x \in B\ |\ wt(x) = n-1\}.$$
It is clear that these subsets are the orbits into which $A$ and $B$
are partitioned under the actions of $S_{p-2}$ and $S_{q-2}$, respectively.

\pr{Lemma 2} For any $x,y\in H$ we have
$$supp(x + y) = supp(x)\cup supp(y) - supp(x)\cap supp(y)$$
is the symmetric difference of the supports of $x$ and $y$.
Furthermore, we have $supp(x)\cap supp(y) = supp(x\cdot y)$, so that
$$wt(x+y) = wt(x) + wt(y) - 2\ wt(x\cdot y).$$
\epr

The following Lemma has an analog for $0 < n_2\leq n_1 < q$.

\pr{Lemma 3} For $0 < m_2\leq m_1 < p$, we have 
$$A_{m_1} + A_{m_2} = \bigcup_{m_3\in M} A_{m_3}$$
where $M = \{m_3\ |\ (m_1,m_2,m_3) \hbox{ is p-admissible}\}$.
\epr

\demo{Proof} To prove that $A_{m_1} + A_{m_2}$ is contained in the
given union, we will prove that
$$\{wt(a_1+a_2)\ |\ a_1\in A_{m_1},\ a_2\in A_{m_2}\} =
\{m_1-m_2+2i\ |\ 0\leq i\leq \min(m_2,p-m_1)-1\}.$$
From Lemma 2 we see that the smallest value of $wt(a_1+a_2)$ occurs when
$supp(a_2)\subseteq supp(a_1)$, in which case $wt(a_1\cdot a_2) = wt(a_2)$ is
maximal. The largest value of $wt(a_1+a_2)$, which occurs when $wt(a_1\cdot
a_2)$ is minimal, is determined in two cases as follows.

{\bf Case 1:} If $wt(a_1)+wt(a_2)\leq p-2$, then the largest value of
$wt(a_1+a_2)$ is $wt(a_1)+wt(a_2)$ when $wt(a_1\cdot a_2) = 0$.

{\bf Case 2:} If $wt(a_1)+wt(a_2)>p-2$ then the largest value of $wt(a_1+a_2)$
is $[(p-2)-wt(a_1)]+[(p-2)-wt(a_2)]=2p-4-(wt(a_1)+wt(a_2))$ when $wt(a_1\cdot
a_2) = wt(a_1)+wt(a_2) - (p-2)$ is minimal.

It is clear that the indicated intermediate values of $wt(a_1+a_2)$ occur as 
$wt(a_1\cdot a_2)$ takes on all intermediate values between its minimum
and maximum. So the union over any proper subset of $M$ would not
contain $A_{m_1} + A_{m_2}$. 

Now, for any $m_3\in M$, and for any $a_3\in A_{m_3}$, we wish to find 
elements in $A_{m_1}$ and $A_{m_2}$ whose sum is $a_3$. From the last paragraph
we see that for some $a_1\in A_{m_1}$ and $a_2\in A_{m_2}$, $a_1 + a_2 =
a'_3\in A_{m_3}$. Let $\sigma\in S_{p-2}$ be a permutation such that
$\sigma(a'_3) = a_3$. Then we have $\sigma(a_1) + \sigma(a_2) = a_3$ with  
$\sigma(a_1)\in A_{m_1}$ and $\sigma(a_2)\in A_{m_2}$. $\hfill\bk$
\enddemo

Let $H_{m,n} = A_m + B_n = \{ a+b\in H\ |\ a \in A_m, b \in B_n\}$.
It is clear that
$$H = \bigcup_{0<m<p}\ \ \bigcup_{0<n<q} H_{m,n}$$
is a disjoint union. Define a surjective map $\tilde{\Phi}: H \to S(p,q)$ by
$\tilde{\Phi}(h)  = (m,n)$ if $h\in H_{m,n}$. Now let $I\approx\boZ_2$ be the
cyclic subgroup of $H$  generated by $(1,1,\cdots ,1)$, and let
$G\approx\boZ_2^{r-1}$ be the quotient group $H/I$. Note that
$H_{m,n}+I=H_{p-m,q-n}+I$, so that we have the induced  surjective map 
$\Phi : G \to V(p,q)$ given by $\Phi(h+I)=[h_{m,n}]$ if $h\in H_{m,n}$.

\pr{Theorem} With notation as above, $\Phi : G \to V(p,q)$ covers the
fusion rules of the $(p,q)$-minimal model.
\epr

\demo{Proof} For any $(m_1,n_1),(m_2,n_2)\in S(p,q)$ we first show that
$$ H_{m_1,n_1} + H_{m_2,n_2} = \bigcup_{(m_3,n_3)\in T} H_{m_3,n_3}$$
where $T = \{(m_3,n_3)\in S(p,q)\ |\ ((m_1,n_1),(m_2,n_2),(m_3,n_3))
\hbox{ is }(p,q)\hbox{-admissible}\}$.
By the $S_3\times S_3$ symmetry of addition in $H = A\times B$, and of
$(p,q)$-admissibility, we may assume that $0 < m_2\leq m_1 < p$ and 
$0 < n_2\leq n_1 < q$. 
The equality is equivalent to showing that
$$A_{m_1} + A_{m_2} = \bigcup_{m_3\in M} A_{m_3} \quad\hbox{ and }\quad
B_{n_1} + B_{n_2} = \bigcup_{n_3\in N} B_{n_3}$$
where $M = \{m_3|(m_1,m_2,m_3) \hbox{ is }p\hbox{-admissible}\}$ and 
$N = \{n_3|(n_1,n_2,n_3) \hbox{ is }q\hbox{-admissible}\}$.
But these follow immediately from Lemma 3 and its analog for $q$.

From this equality in $H$, it is now straightforward to check that $\Phi$ 
covers the fusion rules. $\hfill\bk$ \enddemo

The main result in this paper was presented by A. J. F. on November 18, 1995, at
the Special Session on Quantum Affine Algebras and Related Topics, during the
American  Mathematical Society meeting in Greensboro, NC.

\enddocument
\bye